\documentclass[11pt, a4paper]{amsart}
\usepackage{amsmath,amsthm,amssymb,enumerate,ifpdf}
\usepackage[pdftex]{graphicx}

\numberwithin{equation}{section}
\pdfoutput=1
\newtheorem{thm}{Theorem}[section]

\newtheorem{cor}[thm]{Corollary}
\newtheorem{definition}[thm]{Definition}
\newtheorem{remark}[thm]{Remark}
\newtheorem{example}[thm]{Example}

\newenvironment{df}{\begin{definition}\rm}{\end{definition}}
\newenvironment{rem}{\begin{remark}\rm}{\end{remark}}
\newenvironment{ex}{\begin{example}\rm}{\end{example}}

\newcommand{\acknowledgement}{\noindent {\bfseries Acknowledgements.} }

\title{Braiding surface links which are coverings over the standard torus}

\author {Inasa Nakamura$^*$}
\thanks{$^*$Supported by GCOE, Kyoto University}

\date{}

\address{Research Institute for Mathematical Sciences, Kyoto University, Kyoto 606-8502, 
JAPAN}

 \email{inasa@kurims.kyoto-u.ac.jp}

\begin{document}

\begin{abstract}
We consider a surface link in the 4-space which can be presented by a simple branched covering
over the standard torus, which we call a torus-covering link. 
Torus-covering links include spun $T^2$-knots and turned spun $T^2$-knots. 
In this 
paper we braid a torus-covering link over the standard 2-sphere. 
This gives an upper estimate of the braid index  
of a torus-covering link. 
In particular we show that the turned spun $T^2$-knot of the torus $(2,\,p)$-knot has 
the braid index four. 
\end{abstract}

\maketitle
%

\section{Introduction}
A {\it surface link} is the image of a locally flat embedding of a closed surface into 
the Euclidean 4-space $\mathbb{R}^4$. 
It is known \cite{Kamada1, Kamada3} that any oriented surface link can be presented by the 
closure of a simple surface braid. 
Here, the closure of a surface braid is a surface link of the following form.
Let $S^2$ be a standard 2-sphere in $\mathbb{R}^4$, i.e. the boundary of a standard 3-ball 
in $\mathbb{R}^3 \times \{0\}$. 
The {\it closure of a surface braid} is a surface link embedded in a tubular 
neighborhood $N(S^2)$ of $S^2$ in 
such a way that the projection of it to $S^2$ is a branched covering over $S^2$. 
We identify $N(S^2)$ with $I \times I \times S^2$, where 
$I$ is an interval. 
For a surface link $S$ of such a form,
we consider the singular set of the image of $S$ by the projection to $I \times S^2$, and 
the image of this singular set by the projection to $S^2$ forms a graph on $S^2$. 
An {\it $m$-chart on} $S^2$ is such a graph with certain additional data. 
We can present the original surface link by its $m$-chart on $S^2$ 
(\cite{Kamada2, Kamada3}). 
By removing from $S^2$ a 2-disk which is disjoint with the $m$-chart, 
we obtain an $m$-chart on a 2-disk. 
The resulting $m$-chart on a 2-disk is called a {\it surface link $m$-chart}. 
As we mentioned, any oriented surface link can be presented by the 
closure 
of a simple surface braid. 
Thus it follows that 
any oriented surface link is presented by a surface link $m$-chart. 
In \cite{N}, a \lq\lq torus-covering link" is introduced as 
a new construction of a surface link, 
by considering a standard torus instead of a standard 2-sphere. 
Let $T$ be a standard torus in $\mathbb{R}^4$, i.e. 
the boundary of a standard solid torus in $\mathbb{R}^3 \times \{0\}$. 
A {\it torus-covering link} is a surface link embedded in a tubular 
neighborhood $N(T)$ of $T$ in 
such a way that the projection of it to $T$ is a simple branched covering over $T$. 
For a surface link of such a form,
we can define its {\it $m$-chart on} $T$ in the same way as above. 
By 
cutting $T$ along a meridian and a  
longitude, we obtain an $m$-chart on a 2-disk. 
 We will call the resulting $m$-chart on a 2-disk a {\it torus-covering $m$-chart}. 
A torus-covering link is presented by a torus-covering $m$-chart. 
The aim of this paper is to give a surface link chart description 
of a torus-covering link $S$, 
from the torus-covering $m$-chart which presents $S$. 

Since a torus-covering link is an oriented surface link, it can be presented by 
a surface link $m$-chart. 
We give a surface link chart description 
of a torus-covering link $S$, 
from the torus-covering $m$-chart which presents $S$ (Theorem \ref{theorem3-1}). 
In other words, we braid $S$ over a standard 2-sphere. 
We deform $S$ to the form of 
the closure of a simple surface braid, using the motion picture method. 
The {\it braid index} of an oriented surface link $F$ is the minimum degree of 
simple surface braids whose closures in $\mathbb{R}^4$ are equivalent to $F$. 
 The resulting surface link chart of Theorem \ref{theorem3-1} is a $2m$-chart; 
thus we can see that $2m$ is an upper 
estimate of the braid index of the torus-covering link $S$ (Corollary \ref{c0614-7}). 
In particular, we show that 
the turned spun $T^2$-knot of the torus $(2,p)$-knot has 
the braid index four. 

The paper is organized as follows. 
In Section \ref{s1}, we review the chart description of a simple 
braided surface, and using these terms we give the definition of 
a torus-covering link. 
In Section \ref{braid-index}, we give the statement of Theorem \ref{theorem3-1}. 
In Section \ref{pf}, we prove Theorem \ref{theorem3-1} using the motion picture method. 
In Section \ref{0616-1}, we give an example; we draw the surface link chart presenting 
the turned spun $T^2$-knot of a trefoil. 
%
\section{Torus-covering links} \label{s1}
In this section we give the definition of a torus-covering link. 
In Section \ref{0810-1} we review 
a simple braided surface and its chart description, and in Section \ref{t-c-link} we 
give the definition of 
a torus-covering link. 

\subsection{A braided surface and its chart description} \label{0810-1}
A braided surface was defined in \cite{Rudolph, Kamada3}. 
A surface braid is a braided surface with some boundary condition, and 
a notion of an $m$-chart was introduced 
\cite{Kamada92, Kamada3} to present 
a simple surface braid. 
Equivalent simple surface braids have distinct chart presentations. 
The notion of C-move equivalence between two $m$-charts was introduced 
\cite{Kamada92, Kamada2, Kamada3} to 
give the equivalence class of the chart which represents 
the equivalence class of a simple surface braid. 
The notion of an $m$-chart can be easily extended to an $m$-chart presenting 
a simple braided surface. 
In this subsection we review a braided surface, and extend 
the notion of a chart description to a simple braided surface. 
   
\begin{df}
A compact and oriented 2-manifold $S$ embedded in a bidisk $D_1 \times D_2$ properly 
and locally flatly is called a {\it braided surface} 
of degree $m$
if $S$ satisfies the following conditions:
\begin{enumerate}[(i)]
 \item
 $p_2 | _{S} \,:\, S \rightarrow D_2$ is a branched covering map of 
degree $m$, 
\item 
 $\partial S$ is a closed $m$-braid in $D_1 \times \partial D_2$, 
where $D_1$, $D_2$ are 2-disks, and $p_2 \,:\, D_1 \times D_2
\rightarrow D_2$ is 
the projection to the second factor.
\end{enumerate} 
Two braided surfaces are {\it equivalent} if there is a fiber-preserving ambient isotopy 
of $D_1 \times D_2$ rel $D_1 \times \partial D_2$ which carries one to the other. 
A braided surface $S$ is called \textit{simple} if $\#(S \cap p_2^{-1}(x))=m-1$ or $m$ for each $x \in D_2$. 
A braided surface $S$ is called a \textit{surface braid} if $\partial S$ is the 
trivial closed braid. 
A surface braid $Q_m \times D_2$ is called {\it trivial}, 
where $Q_m$ is a set of $m$ interior points of $D_1$. 
\end{df}

 When a simple braided surface $S$ is given, we obtain a graph on $D_2$, as follows. 
Identify $D_1$ with $I \times I$, where $I=[0,1]$. 
Consider the singular set $\mathrm{Sing}(p_1(S))$ of the image of $S$ by 
the projection $p_1$ to $I \times D_2$. 
Perturbing $S$ if necessary, 
we can assume that 
$\mathrm{Sing}(p_1(S))$ consists of 
double point curves, 
triple points, and branch points. 
Moreover we can assume that 
the singular set of 
the image of $\mathrm{Sing}(p_1(S))$ 
by the projection to 
$D_2$ consists of a finite number of double points such that the preimages 
belong to double point curves of $\mathrm{Sing}(p_1(S))$. 
Thus 
the image of $\mathrm{Sing}(p_1(S))$ by the projection to 
$D_2$ forms a finite graph $\Gamma$ on $D_2$ such that 
the degree of its vertex is either $1$, $4$ or $6$. 
An edge of $\Gamma$ corresponds 
to a double point curve, and a vertex of degree $1$ (resp. $6$) 
corresponds to a branch point (resp. triple point). 

For such a graph $\Gamma$ obtained from a simple braided surface $S$, 
we give orientations and labels to the edges of $\Gamma$, as follows. 
Let us consider a path $\rho$ in $D_2$ such that 
$\rho \cap \Gamma$ is a point $P$ of an edge $e$ of $\Gamma$. 
Then $S \cap p_{2}^{-1} (\rho)$ is a classical $m$-braid with one crossing 
in $p_2^{-1}(\rho)$ 
such that $P$ corresponds to the crossing of the $m$-braid. 
Let $\sigma_1, \sigma_2, \ldots, \sigma_{m-1}$ 
be the standard generators of the $m$-braid group $B_m$. 
Let $\sigma_{i}^{\epsilon}$ 
($i \in \{1,2,\ldots, m-1\}$, 
$\epsilon \in \{+1, -1\}$) be 
the presentation of $S \cap p_{2}^{-1}(\rho)$. 
Then label the edge $e$ by $i$, and moreover give $e$ an orientation such that 
the normal vector of $\rho$ corresponds (resp. 
does not correspond) 
to the orientation of $e$ if $\epsilon=+1$ (resp. $-1$). 
We call such an oriented and labeled graph an {\it $m$-chart of $S$}. 
   \\
 
 In general, we define an $m$-chart on $D_2$ as follows. 

\begin{figure}
 \includegraphics*{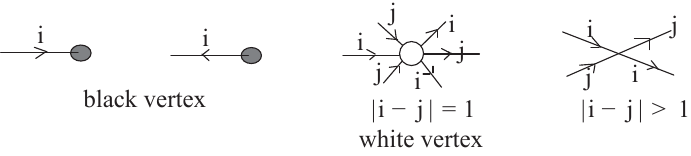}
 \caption{Vertices in a chart.}
\label{Fig1-1}
 \end{figure}

\begin{df} \label{0819-1}
Let $m$ be a positive integer. 
A finite graph $\Gamma$ on a 2-disk $D_2$ is called an {\it $m$-chart} if 
it satisfies 
the following conditions:

\begin{enumerate}[(i)]
\item 
$\Gamma \cap \partial D_2$ consists of a finite number of vertices of degree $1$. 
\item Every edge is oriented and labeled by an element of 
       $\{1,2, \ldots, m-1\}$. 
 \item Every vertex has degree $1$, $4$, or $6$.
 \item  The adjacent edges around each vertex in $\mathrm{Int}(D_2)$ 
are oriented and labeled as shown in 
Figure \ref{Fig1-1}, 
where we depict a vertex of degree 1 by a black vertex, and a vertex of degree 
6 by a white vertex. 
 \end{enumerate}
\end{df}
 A black vertex (resp. white vertex) of an $m$-chart corresponds to a branch point 
 (resp. triple point) of the simple braided surface 
 presented by the $m$-chart. 
An $m$-chart presents a simple braided surface. In particular, an $m$-chart $\Gamma$ such that $\Gamma \cap \partial D_2=\emptyset$ presents a simple surface braid. 
\\

When an $m$-chart $\Gamma$ on $D_2$ is given, 
 we can reconstruct a simple braided surface $S$ over $D_2$ 
as follows. 
Let $N(\Gamma)$ be a neighborhood of $\Gamma$ in $D_2$. 
Let us consider a trivial braided surface $S=Q_m \times (D_2 -N(\Gamma))$ 
over $D_2 -N(\Gamma)$, where $Q_m$ is a set of $m$ interior points of $D_1$. 
We extend $S$ over a neighborhood of each edge as follows. 
Identify a neighborhood of an edge $e$ with $I \times I$ such that 
 $e$ is identified with $\{1/2\} \times I$. 
Let $i$ be the label attached to $e$, and let $\epsilon=+1$ (resp. $-1$) 
if the orientation of $e$ corresponds (resp. does not correspond) 
to the orientation of $\{0\} \times I$. 
Then let the braided surface $S$ over the neighborhood of $e$ be 
the braided surface which has a presentation 
$\sigma_i^\epsilon \times I$ and 
the image of the double point curve of $p_1(S)$ by the projection to $D_2$ 
is $e$. 
Since $\Gamma$ is as in Figure \ref{Fig1-1} around each vertex, 
$S$ can be extended naturally over a neighborhood 
of each vertex. See \cite{Carter-Saito, Kamada92-2, Kamada3} 
for more details. 
Thus we can construct a simple braided surface $S$ over $D_2$ 
such that the original $m$-chart is an $m$-chart of $S$. 
\\

\begin{figure}
 \includegraphics*{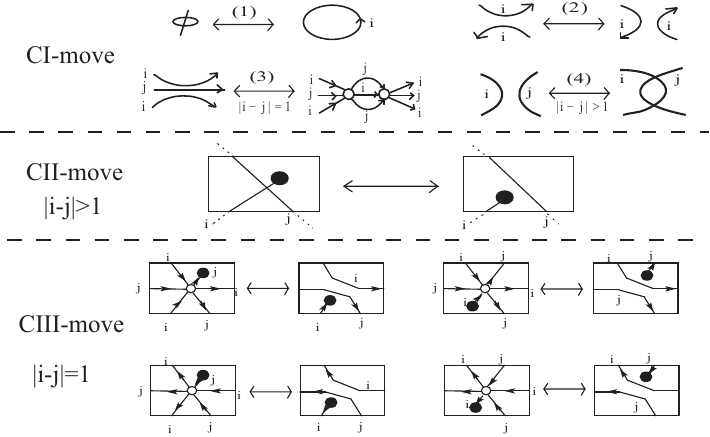}
\caption{CI, CII, CIII-moves. For CI-moves, we give only several examples.} 
\label{cmove}
 \end{figure}

Two $m$-charts on $D_2$ are {\it C-move equivalent} if 
they are related by a finite sequence of 
ambient isotopies of $D_2$ and CI, CII, CIII-moves shown in Figure \ref{cmove}; 
see \cite{Kamada3} for the complete set of CI-moves. 
 It is shown as a minor modification of \cite{Kamada92, Kamada2, Kamada3} that 
two simple braided surfaces of degree $m$ are equivalent if and only if 
$m$-charts of them are C-move equivalent. 
 \\

The boundary of a surface braid $S$ consists of trivial closed $m$-braid. Consider a natural embedding of $D_1 \times D_2$ in $\mathbb{R}^4$, and paste $m$ disks to $S$ to obtain an embedding of a closed surface in $\mathbb{R}^4$. The resulting surface is called the {\it closure} of $S$ (see Section \ref{pf}). 
It is known \cite{Kamada1, Kamada3} that any oriented surface link is presented by the closure of a simple surface braid; thus it is presented by an $m$-chart $\Gamma$ on a 2-disk $D_2$ such that $\Gamma \cap D_2 =\emptyset$. 
We call such an $m$-chart which presents a surface link a {\it surface link $m$-chart}. 

%
\subsection{Torus-covering links} \label{t-c-link}
 
Now we give the definition of a torus-covering link in $\mathbb{R}^4$ (see \cite[Definition 2.1]{N}). 
Let $T$ be a standard torus in $\mathbb{R}^4$, that is, 
the boundary of an unknotted solid torus in a 3-space in $\mathbb{R}^4$. 
Let us consider a tubular neighborhood $N(T)$ of $T$, and identify $N(T)$ with 
$D^2 \times S^1 \times S^1$, 
where 
$D^2$ is a 2-disk, and $S^1$ is a circle. 
The first $S^1$ corresponds to the meridian, and the second $S^1$ corresponds to the 
longitude of $T$. 
Let us identify $S^1$ with $I/\sim$, where $I=[0,\,1]$ and $0 \sim 1$. 
For a manifold $S$ in $N(T)$, let us denote by $S \cap (D^2 \times I \times I)$ 
the manifold in $D^2 \times I \times I$ obtained from $S$ by cutting it 
at $D^2 \times S^1 \times \{0\}$ and $D^2 \times \{0\} \times S^1$. 
\begin{df} \label{2-1} 
  A {\it torus-covering link} is a surface link $S$ in $\mathbb{R}^4$ 
such that 
 $S \subset N(T)$ and moreover 
$S \cap (D^2 \times I \times I)$ 
is a simple braided surface. 
\end{df}
By definition, a torus-covering link $S$ is presented by an $m$-chart 
$\Gamma$ on $I \times I$ with $\Gamma \cap (I \times \{0\})=\Gamma \cap 
(I \times \{1\})$ and $\Gamma \cap (\{0\} \times I)=\Gamma \cap 
(\{1\} \times I)$. 
Let us call $\Gamma$ on $I \times I$ a {\it torus-covering $m$-chart}. 
 
  As we mentioned, 
 for two $m$-charts, 
their presenting braided surfaces 
 are equivalent 
if the $m$-charts are C-move equivalent. 
 Hence it follows that for two torus-covering $m$-charts, 
 their presenting torus-covering links are equivalent 
if the torus-covering $m$-charts are C-move equivalent. 
Since each component of a torus-covering link is a branched cover over a torus $T$, 
 each component of a torus-covering link 
 is of genus at least one. 
 \\

A $T^2$-link is a surface link whose each component is of genus one. 
  As known $T^2$-links constructed from classical links, 
 there are spun $T^2$-links and turned spun $T^2$-links, which are constructed as follows. 
Consider a 3-ball $B^3$ and a natural embedding of $B^3 \times S^1$ in $\mathbb{R}^4$. 
The space $B^3 \times S^1$ can be regarded as constructed by spinning $B^3$ along the circle $S^1$. 
  A {\it spun $T^2$-link} of a classical link $L$ 
is constructed by spinning $L \subset B^3$ (\cite{Livingston, Boyle88, 
Boyle}). 
 A {\it turned spun $T^2$-link} of $L$ is constructed by 
turning $L \subset B^3$ once while spinning (\cite{Boyle}). 
A spun $T^2$-link or a turned spun $T^2$-link 
of any classical closed braid is presented by a torus-covering 
link (see \cite[Propositions 2.10 and 2.11]{N}). 
   
  \begin{figure}
  \includegraphics*{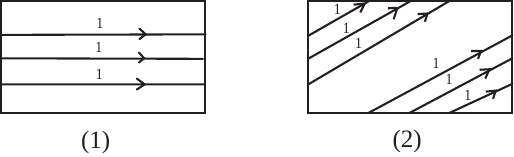}
\caption{Examples of $2$-charts on $T$, or torus-covering $2$-charts on a 2-disk.}
\label{Fig3-10}
\end{figure}
 \begin{ex} \label{0228-1}
 \begin{enumerate}[(1)]
 \item
  Let $\Gamma$ be a torus-covering 2-chart 
 as in Figure \ref{Fig3-10} (1). 
 Then the torus-covering knot 
 presented by $\Gamma$ is the spun $T^2$-knot of 
 a right-handed trefoil.
 \item
  Let $\Gamma$ be 
a torus-covering 2-chart as in Figure \ref{Fig3-10} (2). 
 Then the torus-covering knot 
 presented by $\Gamma$ is the turned spun $T^2$-knot of 
 a right-handed trefoil.
 \end{enumerate}
\end{ex} 

   \section{Braiding a torus-covering link over the 2-sphere} \label{braid-index}
   In this section we give the statement of Theorem \ref{theorem3-1}, 
using the terms of $m$-charts. 
An $m$-chart is presented by a simple braided surface. 
A simple braided surface is presented by 
a motion picture consisting of 
isotopic transformations and hyperbolic transformations (see \cite[Section 9.1]{Kamada3}). 
 
For a subset $A_1$ of $\mathbb{R}^3$ and a subset $A_2$ of $\mathbb{R}$, 
    we denote by $A_1 A_2$ (or by $A_1 \times A_2)$ the subset $\{(x,t) \in 
    \mathbb{R}^3 \times \mathbb{R}\,|\,x \in A_1,\, t \in A_2 \}$ of 
    $\mathbb{R}^3 \times \mathbb{R}=\mathbb{R}^4$. 
    In particular, when $A_2=\{t\}$, we denote $A_1 \{t\}$ by $A_1[t]$. 
    For a subset $X$ of $\mathbb{R}^4=\mathbb{R}^3 \times \mathbb{R}$, 
    a {\it motion picture} of $X$ is a one-parameter family 
    $\{ \pi(X \cap \mathbb{R}^3[t]) \}_{t \in \mathbb{R}}$, 
    where 
    $\pi \,:\, \mathbb{R}^4\rightarrow \mathbb{R}^3$, $\pi(x,y,z,t)=(x,y,z)$ 
    is the projection. 
  
   Let $\{h_t\}_{t \in [0,1]}$ be an ambient isotopy of $\mathbb{R}^3$. For a classical link $L$, 
   we have an isotopy (a one-parameter family) $\{h_t(L)\}$ of classical links. 
   We say that $h_1(L)$ is obtained from $L$ by an {\it isotopic transformation}, 
   and use the notation that $L \rightarrow h_1(L)$ is an isotopic transformation.

   Let $L$ be a classical link in $\mathbb{R}^3$. A 2-disk $B$ in $\mathbb{R}^3$ is called a 
   {\it band} attaching to $L$ if $L \cap B$ is a pair of disjoint arcs in $\partial B$. 
   A {\it band set} attaching to $L$ is a union $\mathcal{B}=B_1 \cup B_2 \cup \cdots \cup B_m$ 
   of mutually disjoint bands $B_1,B_2,\cdots, B_m$ attaching to $L$. 
For a subset $X$ of a space, 
   let us denote by $\mathrm{Cl}(X)$ the closure of $X$. 
   Define a link $h(L;\mathcal{B})$ by 
   \[
   h(L;\mathcal{B})=\mathrm{Cl}\big( (L \cup \partial \mathcal{B})-(L \cap \mathcal{B}) \big). 
   \]
    We say that the link $h(L;\mathcal{B})$ is obtained from $L$ by 
    a {\it hyperbolic transformation} along $\mathcal{B}$, 
    and use the notation that $L \rightarrow h(L;\mathcal{B})$ is a 
    hyperbolic transformation (see \cite[Section 9.1]{Kamada3}).

    For a classical $m$-braid $c$, let 
$\iota^l_k(c)$ be the $(m+k+l)$-braid obtained from $c$ by adding $k$ (resp. $l$) trivial 
strings before (resp. after) $c$, and put 
\begin{eqnarray*}
 & & \Pi ^m_i=\sigma_{m+1} \sigma_{m+2} \cdots \sigma_{m+i},\ 
 \Pi^{\prime \,m}_i=\sigma_{m-1} \sigma_{m-2} 
\cdots \sigma_{m-i},
\\
 & & \Delta_m = \Pi^m_{m-1} \Pi^m_{m-2} \cdots \Pi^m_1, \ 
\Delta_m^\prime = \Pi^{\prime \,m}_{m-1} \Pi^{\prime \,m}_{m-2}
 \cdots \Pi^{\prime \,m}_1,
\\
 & & \Theta_m= \sigma_m \cdot \Pi^{\prime \,m}_{m-1}  \cdot
\Pi^m_{m-1} 
 \cdot \sigma_m \cdot \Pi^{\prime \,m}_{m-2}  \cdot \Pi^m_{m-2}
 \cdots \sigma_m \cdot \Pi^{\prime \,m}_1 
 \cdot \Pi^m_1 \cdot \sigma_m. 
 \end{eqnarray*}
 
\begin{figure}
  \includegraphics*{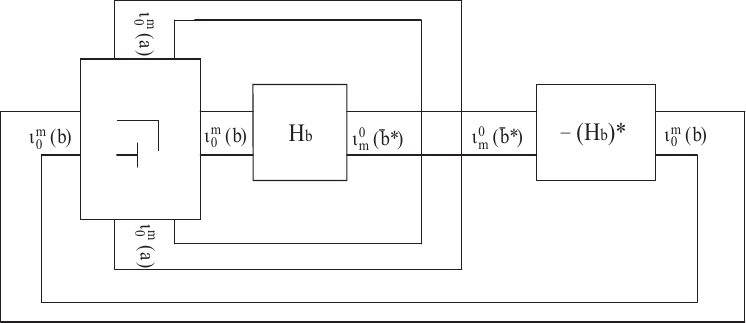}
\caption{The surface link $2m$-chart $\Gamma_S$.}
\label{0614-4}
\end{figure}

\begin{figure}
  \includegraphics*{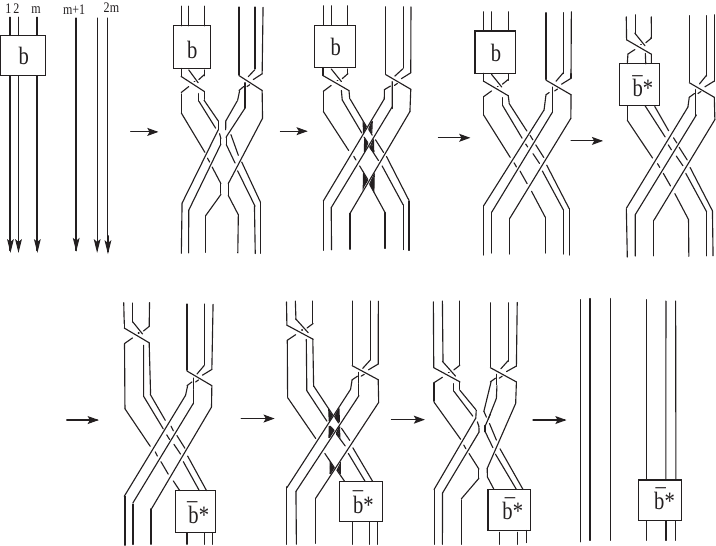}
\caption{The motion picture of the 1-handle braided surface (of degree $2m$)
associated with an $m$-braid $b$.}
\label{Fig3-2}
\end{figure}

\begin{figure}
  \includegraphics*{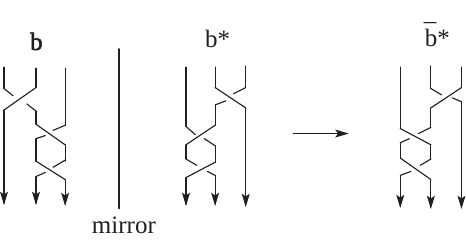}
\caption{The $m$-braid $\bar{b}^*$.}
\label{Fig3-3}
\end{figure}

  \begin{rem}
   Let $\Delta$ be Garside's $\Delta$ (a half twist, see \cite{Garside}) 
  for the $m$-braid group $B_m$. Then $\iota^0_m(\Delta)=\Delta_m$. 
\end{rem}
\begin {thm} \label{theorem3-1}
Let $\Gamma_T$ be a torus-covering $m$-chart. 
Let $a$ (resp. $b$) be a classical $m$-braid presented 
by $\Gamma_T \cap (I \times \{0\})$ (resp. $\Gamma_T \cap (\{0\} \times I)$). 
Then the torus-covering link presented by $\Gamma_T$ 
 is presented by a surface link $2m$-chart $\Gamma_S$ 
 as in Figure \ref{0614-4}.  
 Here 
 $H_b$ is a $2m$-chart 
 presenting the simple braided surface whose motion picture is as follows:
 \begin{eqnarray*}
 \iota^m_0 (b)
 &\longrightarrow &
  \iota^m_0 (b) \cdot (\Delta^{\prime}_m)^{-1} \cdot 
 \Delta_m^{-1} \cdot \Delta^{\prime}_m \cdot
 \Delta_m 
 \dot{\longrightarrow} 
 \iota^m_0(b) \cdot (\Delta^{\prime}_m)^{-1} \cdot 
\Delta^{-1}_m \cdot 
 \Theta_m \\
 &\longrightarrow & 
 (\Delta^{\prime}_m)^{-1} \cdot
\Delta^{-1}_m \cdot \iota^m_0 (\bar{b}^*) \cdot 
 \Theta_m 
 \longrightarrow 
 (\Delta^{\prime}_m)^{-1} \cdot 
 \Delta^{-1}_m \cdot
 \Theta_m \cdot \iota^0_m (\bar{b}^*) \\
 &\dot{\longrightarrow } &
(\Delta^{\prime}_m)^{-1} \cdot 
 \Delta^{-1}_m \cdot 
   \Delta^{\prime}_m \cdot 
\Delta_m \cdot \iota^0_m (\bar{b}^*) 
  \longrightarrow 
  \iota^0_m (\bar{b}^*), 
  \end{eqnarray*}
 where $\longrightarrow$ is an isotopic transformation 
 and $\dot{\longrightarrow}$ is a hyperbolic transformation
  along bands corresponding to the $m\ \sigma_m$'s (see Figure \ref{Fig3-2}), 
 and $-(H_b)^*$
 is the orientation-reversed mirror image of $H_b$,
 and
 $\bar{b}^*$ is the $m$-braid obtained from the classical $m$-braid $b$
 by taking its mirror image and reversing all the crossings (see Figure \ref{Fig3-3}). 
 \end{thm}
 \begin{df}
 Let us call $H_b$ the \textit{1-handle $2m$-chart} associated with an $m$-braid $b$, and its 
 presenting braided surface the 
 \textit{1-handle braided surface}. 
 \end{df}

The surface link $2m$-chart $\Gamma_S$ is well-defined, for 
the edges presenting $\iota^m_0(a)$ 
have labels at most $m-1$ and the edges presenting $\iota^0_m(\bar{b}^*)$ 
have labels at least $m+1$. 
Note that the 1-handle $2m$-chart $H_b$ 
has $2m$ black vertices.

\section{Proof of Theorem \ref{theorem3-1}} \label{pf}
Before proving Theorem \ref{theorem3-1}, we review 
the Cellular Move Lemma and 
the definition of the closure of a surface braid. 

 Let $F$ be a surface link and let $E^3$ be a 3-ball in $\mathbb{R}^4$. 
    Suppose that $E^3 \cap F$ is a 2-disk in $\partial E^3$. 
    Replacing $E^3 \cap F$ by the 2-disk $\mathrm{Cl}\big(\partial E^3 -(E^3 \cap F) \big)$, 
    we obtain a new surface link $F^\prime$ from $F$. Such a replacement is called a {\it cellular move} 
    along $E^3$. If $F$ is oriented, then we assume that $F^\prime$ is oriented by an orientation 
    induced from $\mathrm{Cl}\big(F -(E^3 \cap F) \big)$. 
    It is known that $F$ and $F^\prime$ are equivalent 
({\it Cellular Move Lemma}, see \cite[Lemma 6.6]{Kamada3} or \cite{Rourke-Sanderson}). 
 
     Assume that $L$ is a trivial classical link, and take $t \in \mathbb{R}$. 
    A {\it disk system} in $\mathbb{R}^3(-\infty, t]$ (resp. $\mathbb{R}^3[t, +\infty)$) 
    with boundary $L[t]$ is a union 
    of mutually disjoint 2-disks embedded in 
    $\mathbb{R}^3(-\infty, t]$ (resp. $\mathbb{R}^3[t, +\infty)$) properly and locally flatly 
such that the boundary is 
    $L[t]$.  
    A disk system $\tilde{\mathcal{D}}$ in $\mathbb{R}^3(-\infty, t]$ 
    (resp. $\mathbb{R}^3 [t, +\infty)$) is {\it trivial} if there is no 
    critical point except a single 
    minimal point or minimal disk (resp. a maximal point or maximal disk) 
    on each disk component of $\tilde{\mathcal{D}}$ 
    (see \cite[Section 8.5]{Kamada3}). 
    
    Let $t_1, t_2 \in \mathbb{R}$ such that $t_1<t_2$. 
    Let $F$ be a properly embedded surface in $\mathbb{R}^3[t_1, t_2]$ and let 
    $L_-$ and $L_+$ be links in $\mathbb{R}^3$ such that $L_-$ (resp. $L_+$) 
    appears as the boundary of $F$ 
    in $\mathbb{R}^3[t_1]$ (resp. $\mathbb{R}^3[t_2]$). 
    When $L_-$ and $L_+$ are trivial links, we can consider a closed surface $\hat{F}$ 
    such that 
    \[
    \hat{F}=\tilde{\mathcal{D}}_- \cup F \cup \tilde{\mathcal{D}}_+, 
    \]
    where $\tilde{\mathcal{D}}_-$ (resp. $\tilde{\mathcal{D}}_+$) is a trivial disk system in 
    $\mathbb{R}^3(-\infty, t_1]$ (resp. $\mathbb{R}^3[t_2, +\infty)$) with boundary 
    $L_-[t_1]$ (resp. $L_+[t_2]$). 
    We call $\hat{F}$ the {\it closure} of $F$. 
    It is known \cite[Proposition 9.11]{Kamada3} that for two closures of $F$, 
they are ambient isotopic in $\mathbb{R}^4$ rel 
    $\mathbb{R}^3[t_1, t_2]$. 
    Let $S$ be a surface braid in $D^2 \times I \times I^\prime$, 
    where $D^2$ is a 2-disk, $I=[0,1]$ and $I^\prime$ is an arbitrary interval. 
    Let us consider 
 a map 
    \begin{equation} \label{a}
    A \,:\, D^2 \times I \times I^\prime \rightarrow  
D^2 \times I \times I^\prime/(x, 0, u) \sim (x, 1, u) 
\subset \mathbb{R}^3 \times I^\prime
    \end{equation}
    ($x \in D^2$ and $u \in I^\prime$), 
where $A(D^2 \times I \times I^\prime)=D^2 \times S^1 \times I^\prime 
 \subset \mathbb{R}^3 \times I^\prime$
is a natural embedding, i.e. 
    $D^2 \times S^1$ is naturally embedded in $\mathbb{R}^3$. 
    The image 
    $A(S)$ is embedded in $\mathbb{R}^3 \times I^\prime$ properly and locally flatly. 
    The {\it closure} of the surface braid $S$ is the closure of $A(S)$. 
   \\

\noindent
{\it Proof of Theorem \ref{theorem3-1}}. 
Let $S$ be the torus-covering link presented by a torus-covering $m$-chart $\Gamma_T$, and 
let $\pi \,:\, \mathbb{R}^4\rightarrow \mathbb{R}^3$, $\pi(x,y,z,t)=(x,y,z)$ 
    be the projection. 
Let us consider 
the motion picture $\{S_t\}_{t \in \mathbb{R}}$ of $S$, 
where $S_t=\pi(S \cap \mathbb{R}^3 [t])$. 
 We will draw the figures of $S_t$ as 
the diagrams projected on the xy-plane. 
We will use the same notation $S$ after each equivalent deformation. 
For a braid $c$, let us denote by $\mathrm{cl}(c)$ the closed braid 
associated with $c$.

\begin{figure} 
     \includegraphics*{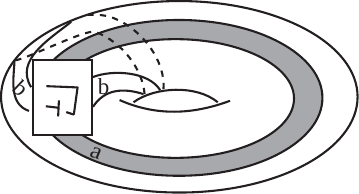}
  \caption{The torus-covering link $S$.}
  \label{B1-10}
\end{figure} 

\begin{figure} 
     \includegraphics*{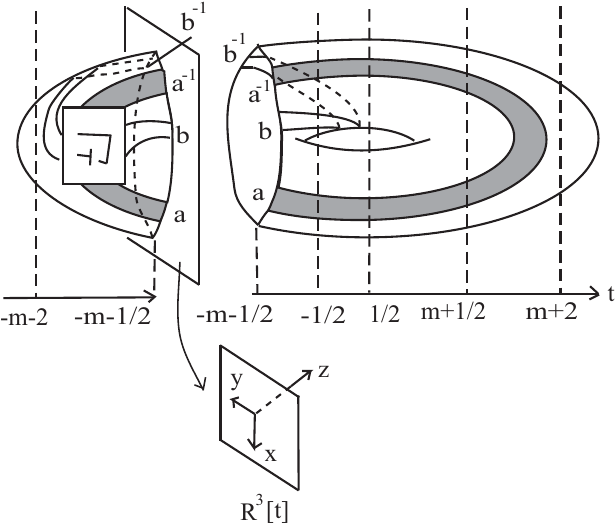}
  \caption{The torus-covering link $S$.}
  \label{B1-11}
\end{figure} 

\begin{figure} 
     \includegraphics*{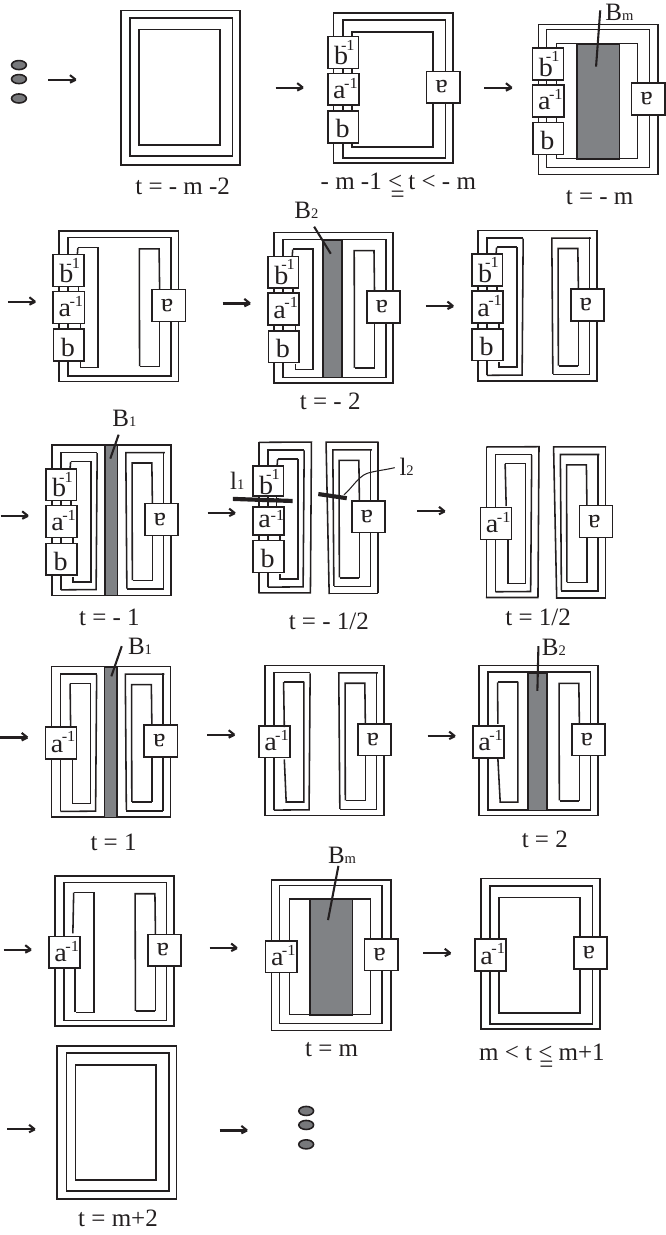}
  \caption{Motion picture of $S$ (1).}
  \label{B1-3}
\end{figure}  

\begin{figure} 
     \includegraphics*{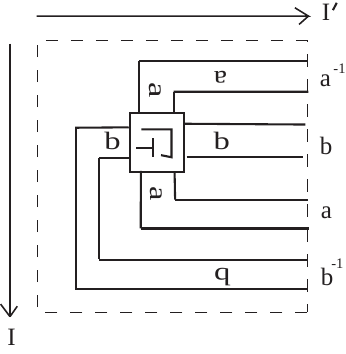}
  \caption{The $m$-chart $\Gamma_T^\prime$.}
  \label{B4-2}
\end{figure} 
{\bf Step 1}. \ 
 Let us consider 
the standard torus $T$ in the $xyt$-space with 
$\Gamma_T$ on it (see Figure \ref{B1-10}). 
We can assume that $T$ with $\Gamma_T$ is as in Figure \ref{B1-11}. 
Then the motion picture of $S$ is as in Figure \ref{B1-3}, as follows. 

Let $B_1, B_2, \ldots, B_m$ 
be $m$ bands in $\mathbb{R}^3$ 
such that $B_j[-j]$ (resp. $B_j[j]$) is an untwisted band 
attaching to $S_{-j-1/2}[-j] \subset \mathbb{R}^3[-j]$ (resp. $S_{j-1/2}[j] \subset \mathbb{R}^3[j]$) 
as in Figure \ref{B1-3} ($j=1,2,\ldots,m$). 
Let us denote by $e_m$ the trivial $m$-braid. 
 Let $\tilde{\mathcal{D}}_{-,m}$ (resp. $\tilde{\mathcal{D}}_{+, m}$) be 
 a trivial disk system of $m$ components in 
$\mathbb{R}^3(-\infty, t]$ (resp. $\mathbb{R}^3[t, +\infty)$) with boundary 
$\mathrm{cl}(e_m)[t]$, 
where $t \in \mathbb{R}$. 
 Let $\Gamma_T^\prime$ be an $m$-chart on $I \times I^\prime$ as in Figure \ref{B4-2}, 
and 
 let 
 $S(\Gamma_T^\prime)$ be the simple braided surface of degree $m$ in 
$D^2 \times I \times I^\prime$ 
 presented by $\Gamma_T^\prime$. 
 Then let us consider $A(S(\Gamma_T^\prime))$, where $A \,:\, 
 D^2 \times I \times I^\prime \rightarrow \mathbb{R}^3 \times I^\prime$ is the map (\ref{a}), 
and 
 denote this by $A(\Gamma_T)$. 
 Note that $A(\Gamma_T)$ is embedded 
 in $\mathbb{R}^3 \times I^\prime$ properly and locally flatly. 
Let us call $A(\Gamma_T)$ the 
{\it braided surface over an annulus} associated with $\Gamma_T$. 
Let us consider an $m$-chart $\Gamma$. 
 Let $\rho$ be a path 
which intersects with the edges of $\Gamma$ transversely. 
Let $c=\sigma_{i_1}^{\epsilon_{i_1}} \sigma_{i_2}^{\epsilon_{i_2}}\cdots \sigma_{i_\nu}^{\epsilon_{i_\nu}}$ 
be the $m$-braid presented by $\rho \cap \Gamma$. 
Let us consider the orientation-reversed path of $\rho$, 
and denote it by $-\rho$. Then 
 $-\rho \cap \Gamma$ presents 
$\sigma_{i_\nu}^{-\epsilon_{i_\nu}}\sigma_{i_{\nu-1}}^{-\epsilon_{i_{\nu-1}}}\cdots \sigma_{i_1}^{-\epsilon_{i_1}}$, 
which is $c^{-1}$. 
Hence the classical $m$-braid $S(\Gamma_T^\prime) \cap (D^2 \times I \times \{ t\})$ 
is $a^{-1}bab^{-1}$, 
where $\{t\}$ is the point of $\partial I^\prime$ with the greater coordinate. 
 
Let $l_1$ (resp. $l_2$) $\subset \mathbb{R}^3$ be 
a half plane 
indicated at $t=-1/2$ 
   in Figure \ref{B1-3}. 
Let us take the identified corresponding ends of a closed braid in 
$l_1$ (resp. $l_1 \cup l_2$) for $t \in [-m-2,\, -m) \cup (m,\, m+2]$ (resp. $t \in (-1, 1)$). 
Then the motion picture depicted in Figure \ref{B1-3}, which describes $S$ 
presented by $\Gamma_T$, 
is as follows. 
A {\it split union} of two classical links $L_1$ and $L_2$ in $\mathbb{R}^3$ is a 
classical link presented by the union of the copies of $L_1$ and $L_2$ 
such that for a 2-sphere $S^2$ embedded in $\mathbb{R}^3$, 
$L_1$ is inside of $S^2$ and $L_2$ is outside. 
 We use $\cup$ for a split union of classical links. 
\begin{enumerate}
\item
$S \cap \mathbb{R}^3 (-\infty , -m-2]=\tilde{\mathcal{D}}_{-,m}$. 
We have 
$S_{-m-2}=\mathrm{cl}(e_m)$. 

\item 
 $S \cap \mathbb{R}^3 [-m-2, -m-1]=A(\Gamma_T)$ 
such that $S_{-m-1}=\mathrm{cl}(a^{-1}bab^{-1})$. 

\item 
$S_t=S_{-m-1}$, for $t \in [-m-1, \,-m-1/2]$. 
 
 \item
$S_{-j-1/2} \rightarrow S_{-j+1/2}$ is a hyperbolic transformation along 
the band $B_j \subset S_{-j}$, 
where $j=1,2,\ldots,m$. 
We have $S_{-1/2}=\mathrm{cl}(a^{-1} b b^{-1}) \cup \mathrm{cl}(a)$. 

\item
$S_{-1/2} \rightarrow S_{1/2}$ is an isotopic transformation such that 
$S_{1/2}=\mathrm{cl}(a^{-1}) \cup \mathrm{cl}(a)$. 

\item 
$S_{j-1/2} \rightarrow S_{j+1/2}$ is a hyperbolic transformation along 
the band $B_j \subset S_j$, 
where $j=1,2,\ldots,m$. 
We have $S_{m+1/2}=\mathrm{cl}(a^{-1}a)$. 

\item 
$S_t=S_{m+1/2}$ for $t \in [m+1/2, m+1]$. 

\item
$S_{m+1} \rightarrow S_{m+2}$ is an isotopic transformation such that 
$S_{m+2}=\mathrm{cl}(e_m)$. 

\item 
$S \cap \mathbb{R}^3 [m+2, +\infty)=\tilde{\mathcal{D}}_{+,m}$. 
\end{enumerate}
 
\begin{figure} 
     \includegraphics*{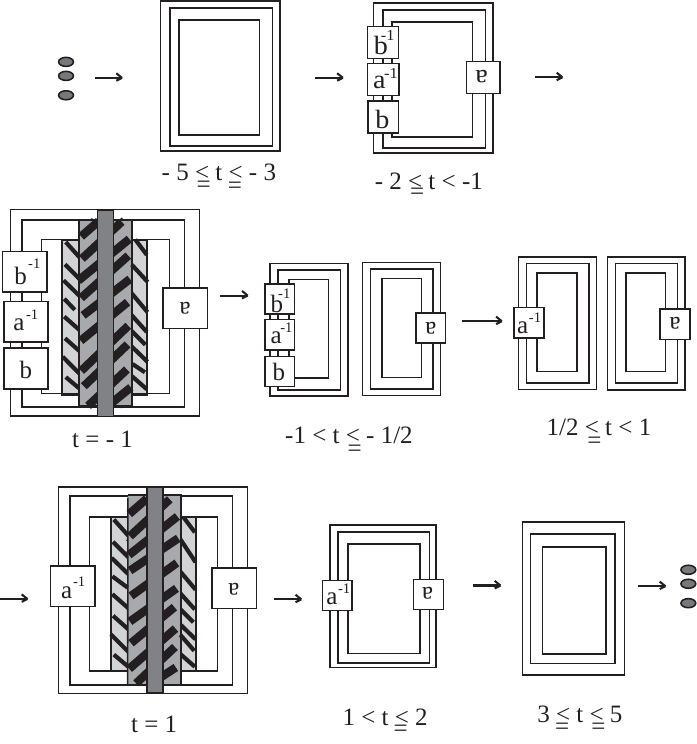}
  \caption{Motion picture of $S$ (2). }
  \label{0614-2}
\end{figure} 

{\bf Step 2}. \ 
Now, let us move bands $B_2[-2], B_3[-3], \ldots, B_m[-m]$ (resp. 
$B_2[2], B_3[3], \ldots, B_m[m]$) into $\mathbb{R}^3 [-1]$ (resp. 
$\mathbb{R}^3 [1]$) 
as follows. 
For subsets $X$ and $Y$ of the $xyz$-space $\mathbb{R}^3$, 
we will say {\it $X$ is over $Y$ with respect to the $z$-axis} 
if for any pair $((x,y,z_1), (x,y,z_2)) \in X \times Y$, $z_1>z_2$ holds. 
 
Now we can assume that $B_{j-1}$  is over $B_{j}$ 
with respect to the $z$-axis, where $j=2,3,\ldots,m$. 
Let us consider 3-balls 
$V_2^\epsilon, V_3^\epsilon, \ldots, V_m^\epsilon$ ($\epsilon=-1, +1$)
in $\mathbb{R}^4$ 
 such that 
\begin{eqnarray*}
\pi(V_j^\epsilon \cap \mathbb{R}^3 [t])=
\begin{cases}
B_j & \mathrm{for}\ t \in [-j, -1] \ \mathrm{if} \ \epsilon=-1, \ \mathrm{or} \\
      & \mathrm{for} \ t \in [1,j] \ \mathrm{if} \ \epsilon=+1, \\
\emptyset & \mathrm{otherwise}, 
\end{cases}
\end{eqnarray*}
where $j=2,3,\ldots,m$. 
Since $B_{j-1}$ is over $B_{j}$ with respect to the $z$-axis, we can see that 
$V_2^\epsilon, V_3^\epsilon, \ldots, V_m^\epsilon$ 
are mutually disjoint and $V_j^\epsilon \cap S$ is the 2-disk 
$B_j[\epsilon j] \subset \mathbb{R}^3[\epsilon j]$, which is in 
$\partial V_j^\epsilon$, 
where $\epsilon=-1, +1$, $j=2,3,\ldots,m$. 
Hence, by cellular moves along the 3-balls $V_2^\epsilon, V_3^\epsilon,\ldots, V_m^\epsilon$, 
we have a band set $\mathcal{B}_1[\epsilon]$ 
in $\mathbb{R}^3 [\epsilon]$ ($\epsilon=-1, +1$), where 
$\mathcal{B}_1=B_1 \cup B_2\cup \cdots \cup B_m$. 
Let us use the same notation $S$ for the new surface link obtained from $S$ by the cellular moves. 
By the Cellular Move Lemma, they are equivalent. 
 Then we have the motion picture of $S$ as follows: 
 \begin{enumerate}
 \item
$S \cap \mathbb{R}^3 (-\infty , -m-2]=\tilde{\mathcal{D}}_{-,m}$. 
We have  
$S_{-m-2}=\mathrm{cl}(e_m)$.
 
\item 
 $S \cap \mathbb{R}^3 [-m-2, -m-1]=A(\Gamma_T)$. 
We have $S_{-m-1}=\mathrm{cl}(a^{-1}bab^{-1})$.

 \item
$S_{-m-1} \rightarrow S_{-1/2}$ is a hyperbolic transformation along 
the band set $\mathcal{B}_1 \subset S_{-1}$. 
We have $S_{-1/2}=\mathrm{cl}( a^{-1} b b^{-1}) \cup \mathrm{cl}(a)$. 

\item
$S_{-1/2} \rightarrow S_{1/2}$ is an isotopic transformation 
such that 
$S_{1/2}=\mathrm{cl}(a^{-1}) \cup \mathrm{cl}(a)$. 

\item 
$S_{1/2} \rightarrow S_{m+1}$ is a hyperbolic transformation along 
the band set $\mathcal{B}_1 \subset S_1$. 
We have $S_{m+1}=\mathrm{cl}(a^{-1}a)$. 

\item 
$S_{m+1} \rightarrow S_{m+2}$ is an isotopic transformation such that 
$S_{m+2}=\mathrm{cl}(e_m)$. 

\item 
$S \cap \mathbb{R}^3 [m+2, +\infty)=\tilde{\mathcal{D}}_{+,m}$. 
\end{enumerate}
 
Then, by an ambient isotopy of $\mathbb{R}^4$ rel 
$\mathbb{R}^3[-2, 2]$, we can deform 
$S$ 
to have the motion picture as in Figure \ref{0614-2}, 
which is as follows: 
 \begin{enumerate}
  \item
$S \cap \mathbb{R}^3 (-\infty , -5]=\tilde{\mathcal{D}}_{-,m}$.  
We have 
$S_{-5}=\mathrm{cl}(e_m)$.
 
 \item
 $S_t=S_{-5}$ for $t \in [-5, -3]$. 
 
\item 
 $S \cap \mathbb{R}^3 [-3, -2]=A(\Gamma_T)$. 
We have $S_{-2}=\mathrm{cl}(a^{-1}bab^{-1})$. 
 
 \item
$S_{-2} \rightarrow S_{-1/2}$ is a hyperbolic transformation along 
the band set $\mathcal{B}_1 \subset S_{-1}$. 
We have $S_{-1/2}=\mathrm{cl}( a^{-1} b b^{-1}) \cup \mathrm{cl}(a)$. 

\item
$S_{-1/2} \rightarrow S_{1/2}$ is an isotopic transformation 
such that 
$S_{1/2}=\mathrm{cl}(a^{-1}) \cup \mathrm{cl}(a)$. 

\item 
$S_{1/2} \rightarrow S_{2}$ is a hyperbolic transformation 
along the band set $\mathcal{B}_1 \subset S_1$. 
We have $S_{2}=\mathrm{cl}(a^{-1}a)$. 

\item 
$S_{2} \rightarrow S_{3}$ is an isotopic transformation such that 
$S_{3}=\mathrm{cl}(e_m)$. 

\item
$S_t=S_3$ for $t \in [3,5]$. 

\item 
$S \cap \mathbb{R}^3 [5, +\infty)=\tilde{\mathcal{D}}_{+,m}$. 
\end{enumerate}
 
\begin{figure} 
     \includegraphics*{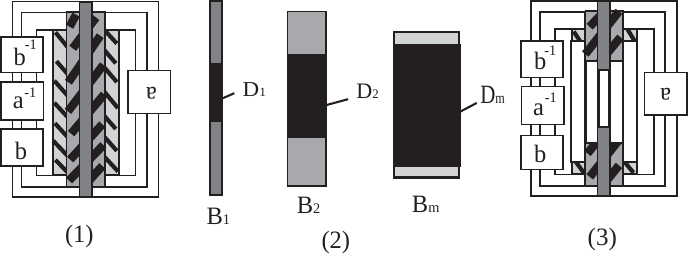}
  \caption{(1) $S_{-1}$ before the cellular moves; 
  (2) 2-disks $D_i$ in $B_i$; (3) $S_{-1}$ after the cellular moves.}
  \label{B1-5}
\end{figure} 

\begin{figure} 
     \includegraphics*{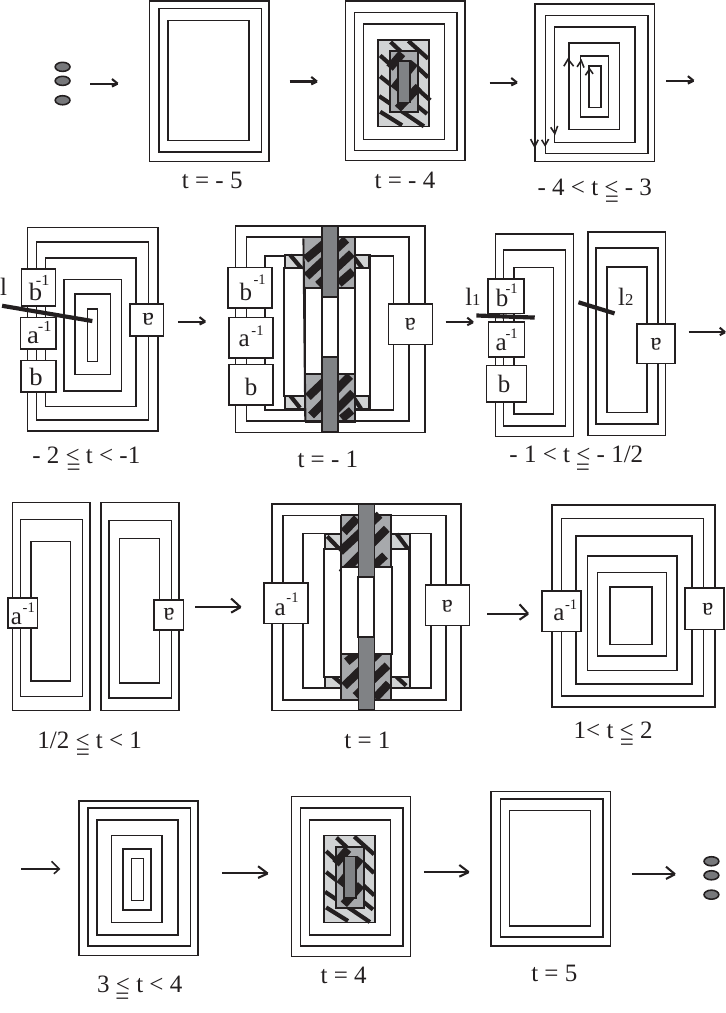}
  \caption{Motion picture of $S$ (3). }
  \label{0614-1}
\end{figure}   

\begin{figure} 
     \includegraphics*{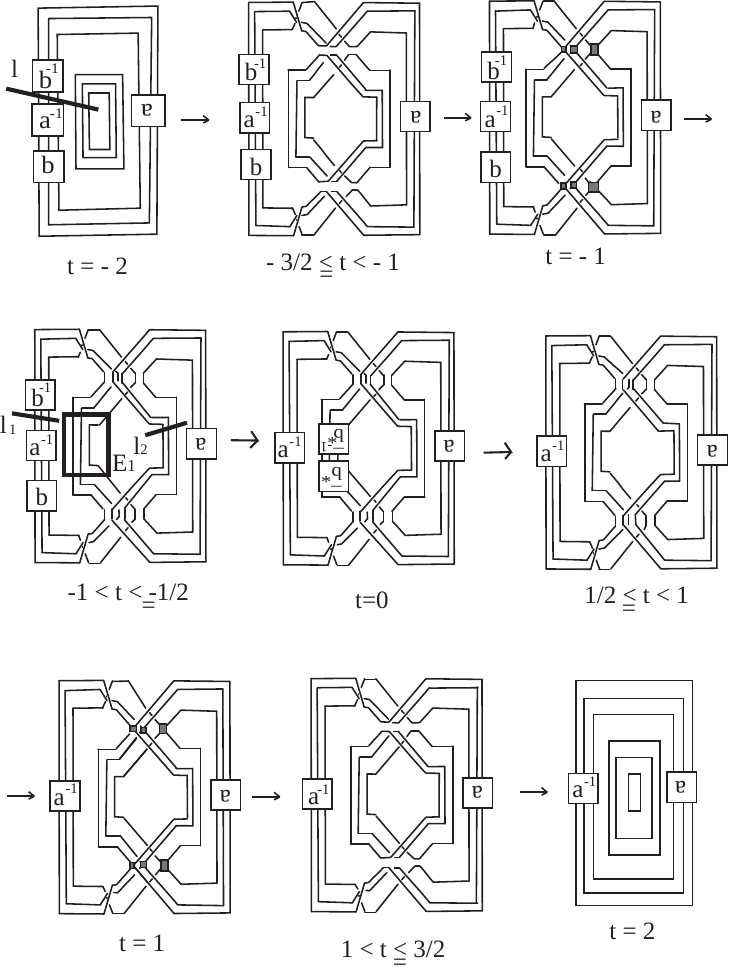}
  \caption{Motion picture of $S \cap \mathbb{R}^3[-2,2]$ (4). }
  \label{0613-1}
\end{figure}  

{\bf Step 3}. \  
We will show that $S$ is equivalent to the surface link whose 
motion picture $\{S_t\}$ is as in Figure \ref{0614-1} for $t\in (-\infty, -2) \cup (2, +\infty)$, 
and otherwise as in Figure \ref{0613-1}.  

Let us take 2-disks $D_1, D_2,\ldots,D_m$ in $\mathbb{R}^3$ 
such that 
$D_j \subset B_j$ 
and $\mathrm{Cl}(B_j-D_j)$ consists of two 2-disks such that each  
contains one of the two disjoint arcs of $S_{-3/2} \cap B_j$; 
see (2) of Figure \ref{B1-5}. 
 Let $W_j^\epsilon$ ($\epsilon=-1, +1$) be a 3-ball in $\mathbb{R}^4$ such that 
 \begin{eqnarray*}
\pi(W_j^\epsilon \cap \mathbb{R}^3 [t])=
\begin{cases}
D_j & \mathrm{for}\ t \in [-4, -1] \ \mathrm{if} \ \epsilon=-1, \ \mathrm{or} \\
    & \mathrm{for}\ t \in [1,4] \ \mathrm{if} \ \epsilon=+1, \\
\emptyset & \mathrm{otherwise}, 
\end{cases}
\end{eqnarray*}
where $j=1,2,\ldots,m$. 
Since $B_1, B_2, \ldots, B_m$ 
are mutually disjoint, 
so are 
$D_1, D_2,\ldots, D_m$; thus 
$W_1^\epsilon, W_2^\epsilon,\ldots, W_m^\epsilon$
 are mutually disjoint. 
Moreover we can see that 
$W_j^\epsilon \cap S$ is the 2-disk $D_j[\epsilon] \subset \mathbb{R}^3[\epsilon]$, 
which is in $\partial W_j^\epsilon$,  
where $\epsilon=-1, +1$, $j=1,2,\ldots,m$. 
By cellular moves along the 3-balls $W_1^\epsilon, W_2^\epsilon,\ldots, W_m^\epsilon$
 ($\epsilon=-1, +1$), 
we have $S$ whose motion picture is as in Figure \ref{0614-1}. 
Figure \ref{B1-5} shows (1) $S_{-1}$ before the cellular moves, 
(2) 2-disks $D_1, D_2,\ldots,D_m$ in $B_1, B_2, \ldots, B_m$, and 
(3) $S_{-1}$ after the cellular moves. 
Put $\mathcal{B}_2=\mathrm{Cl} (B_1-D_1) \cup \mathrm{Cl} (B_2-D_2) \cup \cdots \cup 
\mathrm{Cl} (B_m-D_m)$, which 
is a band set in $\mathbb{R}^3$ consisting of $2m$ bands. 
Then the motion picture of $S$ as in Figure \ref{0614-1} is as follows. 
 Let $-\tilde{\mathcal{D}}_{-,m}$, $-\tilde{\mathcal{D}}_{+,m}$, or $-\mathrm{cl}(e_m)$ be 
 the orientation-reversed image of $\tilde{\mathcal{D}}_{-,m}$, $\tilde{\mathcal{D}}_{+,m}$ or  
$\mathrm{cl}(e_m)$ respectively. 
For an $m$-braid $c$, let us denote by $\iota c$ the braid $\iota_m^0(c)-e_m$, 
and by $-\iota \tilde{\mathcal{D}}_{-,m}$ (resp. $-\iota \tilde{\mathcal{D} }_{+,m}$) 
the trivial disk system of $m$ components in 
$\mathbb{R}^3 (- \infty, t]$
 (resp. $\mathbb{R}^3 [t, +\infty)$) with boundary $-\iota e_m[t]$, 
where $t \in \mathbb{R}$. 
Let $l$ (resp. $l_1$, $l_2$) $\subset \mathbb{R}^3$ be a half plane 
indicated at $t=-2$ (resp. $-1/2$) 
   in Figure \ref{0614-1}. 
Let us take the identified corresponding ends of a closed braid in 
$l$ (resp. $l_1 \cup l_2$) for $t \in [-5,\, -1) \cup (1,\, 5]$ (resp. $t \in (-1, 1)$). 
\begin{enumerate}
\item
$S \cap \mathbb{R}^3 (-\infty , -5]=\tilde{\mathcal{D}}_{-,m}$. 
We have $S_{-5}=\mathrm{cl}(e_m)$. 
\item
$S_t=S_{-5}$ for $t \in [-5, -4)$. 

\item
$S_{-4}=\mathrm{cl}(e_m) \cup -\iota \tilde{\mathcal{D}}_{-,m}$.  

\item
$S_t=
\mathrm{cl}(e_m) \cup -\mathrm{cl}(\iota e_m)$ for $t \in (-4,-3]$. 

\item 
$S \cap \mathbb{R}^3 [-3, -2]=A(\Gamma_T) \cup -\mathrm{cl}(\iota e_m) [-3, -2]$. 
We have $S_{-2}=
\mathrm{cl}(a^{-1}bab^{-1}) \cup -\mathrm{cl}(\iota e_m)$. 

 \item
$S_{-2} \rightarrow S_{-1/2}$ is a hyperbolic transformation along 
the band set $\mathcal{B}_2 \subset S_{-1}$. 
We have $S_{-1/2}=\mathrm{cl}(b^{-1} a^{-1} b) \cup \mathrm{cl}(a)$. 

\item
$S_{-1/2} \rightarrow S_{1/2}$ is an isotopic transformation 
such that 
$S_{1/2}=\mathrm{cl}(a^{-1}) \cup \mathrm{cl}(a)$. 

\item 
$S_{1/2} \rightarrow S_{2}$ is a hyperbolic transformation along 
the band set $\mathcal{B}_2 \subset S_1$. 
We have $S_{2}=\mathrm{cl}(aa^{-1}) \cup -\mathrm{cl}(\iota e_m)$. 

\item 
$S_{2} \rightarrow S_3$ 
is an isotopic transformation such that 
$S_3=\mathrm{cl}(e_m) \cup -\mathrm{cl}(\iota e_m)$. 

\item
$S_t=S_3$ for $t \in (3,4)$. 

\item 
$S_4=
\mathrm{cl}(e_m) \cup -\iota \tilde{\mathcal{D}}_{+,m}$.  

\item 
$S_t=e_m$ for $t \in (4, 5]$.

\item 
$S \cap \mathbb{R}^3 [5, +\infty)=\tilde{\mathcal{D}}_{+,m}$. 
\end{enumerate}

Let $\Delta$ (resp. $\Delta^{\prime}$) be the $m$-braid obtained from 
$\Delta_m$ (resp. $\Delta_m^\prime$) by removing $m$ trivial strings which are 
from the first string to the $m$th string (resp. from the $(m+1)$th string to the $2m$th string), 
 and let 
$\bar{\Delta}$ be the $m$-braid obtained from $\Delta$ by changing each crossing. 
Since $\mathrm{Cl}(B_{j-1}-D_{j-1})$ is over $\mathrm{Cl}(B_{j}-D_{j})$ with respect to the $z$-axis
($j=2,3,\ldots,m$), 
we can take 
an ambient isotopy $\{f_u\}_{u \in [0,1]}$ of $\mathbb{R}^3$ satisfying the following conditions. 
\begin{itemize}
\item
$f_0=\mathrm{id}$. 
\item
$f_u$ is relative the complement of 
a neighborhood of $\cup_{j=1}^m \mathrm{Cl}(B_j-D_j)$, and 
$f_u(\mathrm{cl}(e_m) \cup -\mathrm{cl}(\iota e_m))$ is a split union of 
closed braids for each $u \in [0,1]$, i.e. 
$f_u(S_t)$ is a split union of 
closed braids for each $u \in [0,1]$ and $t \in [-2, 2]-\{-1, 1\}$.
\item 
$f_1 \big( \mathrm{cl}(e_m) \cup -\mathrm{cl}(\iota e_m) \big)=
\mathrm{cl}(\Delta^{\prime -1} \Delta^{\prime} \Delta^{\prime -1} \Delta^{\prime}) \cup 
-\mathrm{cl}(\iota \bar{\Delta} \bar{\Delta}^{-1})
$
such that 
the two untwisted bands $f_1(\mathrm{Cl}(B_j-D_j)) \subset f_1(\mathcal{B}_2)$ 
connect the $j$th string of 
$\mathrm{cl}(\Delta^{\prime -1} \Delta^{\prime} \Delta^{\prime -1} \Delta^{\prime})$ 
with the $j$th string of 
$-\mathrm{cl}(\iota \bar{\Delta} \bar{\Delta}^{-1})$, where $j=1,2,\ldots,m$. 
Note that the $(j-1)$th string of 
$\mathrm{cl}(\Delta^{\prime -1} \Delta^{\prime} \Delta^{\prime -1} \Delta^{\prime})$ 
or $-\mathrm{cl}(\iota \bar{\Delta} \bar{\Delta}^{-1})$ is over 
the $j$th string with respect to the $z$-axis ($j=2,3,\ldots,m$). 
\end{itemize}
Let us take such an ambient isotopy $\{f_u\}$. 

Let $\{h_u\}_{u \in [0,1]}$ be an ambient isotopy of $\mathbb{R}^3$. 
Let us cosider an ambient isotopy $\{H_u\}_{u \in [0,1]}$ of $\mathbb{R}^3[t_1, t_4]$ 
rel $\mathbb{R}^3[t_1] \cup \mathbb{R}^3[t_4]$ 
such that for each $u \in [0,1]$ and $\mathbf{x} \in \mathbb{R}^3$, 
\begin{equation} \label{0528-9}
H_u(\mathbf{x}, t)=\begin{cases}
\big(h_{u(t-t_1)/(t_2-t_1)}(\mathbf{x}), t \big) & \mathrm{if} \ t \in [t_1, t_2] \\
\big(h_u(\mathbf{x}), t \big) & \mathrm{if} \ t \in [t_2, t_3] \\
\big(h_{u(t_4-t)/(t_4-t_3)}(\mathbf{x}), t \big) & \mathrm{if} \ t \in [t_3, t_4], 
\end{cases}
\end{equation}
where $t_1<t_2<t_3<t_4$. 
Then by an ambient isotopy of (\ref{0528-9}) with $h_u=f_u$ and 
$t_1=-2$, $t_2=-3/2$, $t_3=3/2$, $t_4=2$, 
$S$ is equivalent to the surface link whose motion picture $\{S_t\}$ is 
as in Figure \ref{0614-1} 
for $t \in (-\infty, -2) \cup (2, \infty)$, and otherwise 
as follows (see Figure \ref{0613-1}), 
where $\mathcal{B}_3
=f_1(\mathcal{B}_2)$ and 
we take the identified corresponding ends of a closed braid in 
$l$ (resp. $l_1 \cup l_2$) for $t \in [-2,\, -1) \cup (1,\, 2]$ (resp. $t \in (-1, 1)$). 
Here $l$ (resp. $l_1$, $l_2$) 
is a half plane 
indicated at $t=-2$ (resp. $t=-1/2$) in Figure \ref{0613-1}. 
\begin{enumerate}
\item 
$S_{-2} \rightarrow S_{-3/2}$ is an isotopic transformation by $\{f_t\}$. 
We have $S_{-3/2}=\mathrm{cl}(a^{-1}b \Delta^{\prime -1} \Delta^\prime a 
\Delta^{\prime -1} \Delta^\prime b^{-1}) 
\cup -\mathrm{cl}(\iota \bar{\Delta} \bar{\Delta}^{-1})$.  

 \item
$S_{-3/2} \rightarrow S_{-1/2}$ is a hyperbolic transformation along 
the band set $\mathcal{B}_3 \subset S_{-1}$. 
We have $S_{-1/2}=\mathrm{cl}(a^{-1} b \Delta^{\prime -1} \Delta^\prime b^{-1}) 
\cup \mathrm{cl}(a)$.  

\item
$S_{-1/2} \rightarrow S_{1/2}$ is an isotopic transformation such that 
$S_0=\mathrm{cl}(a^{-1} \Delta^{\prime -1} \bar{b}^* \bar{b}^{* -1} \Delta^\prime) 
\cup \mathrm{cl}(a)$ and 
$S_{1/2}=\mathrm{cl}(a^{-1} \Delta^{\prime -1} \Delta^\prime) \cup \mathrm{cl}(a)$.  

\item 
$S_{1/2} \rightarrow S_{3/2}$ is a hyperbolic transformation along the band set 
$\mathcal{B}_3 \subset S_1$. We have $S_{3/2}=\mathrm{cl}(a^{-1} \Delta^{\prime -1} 
\Delta^\prime a \Delta^{\prime -1} \Delta^\prime) 
\cup -\mathrm{cl}(\iota \bar{\Delta} \bar{\Delta}^{-1})$. 

\item 
$S_{3/2} \rightarrow S_{2}$ is an isotopic transformation by $\{f_{1-t} \}$.  
We have $S_2=\mathrm{cl}(a^{-1} a) \cup -\mathrm{cl}(\iota e_m)$. 
\end{enumerate}
  
 \begin{figure} 
     \includegraphics*{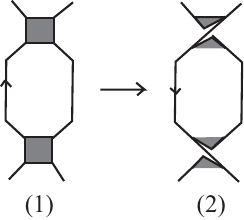}
  \caption{The trivial closed braid, arcs and bands.}
  \label{0528-2}
\end{figure} 

\begin{figure} 
     \includegraphics*{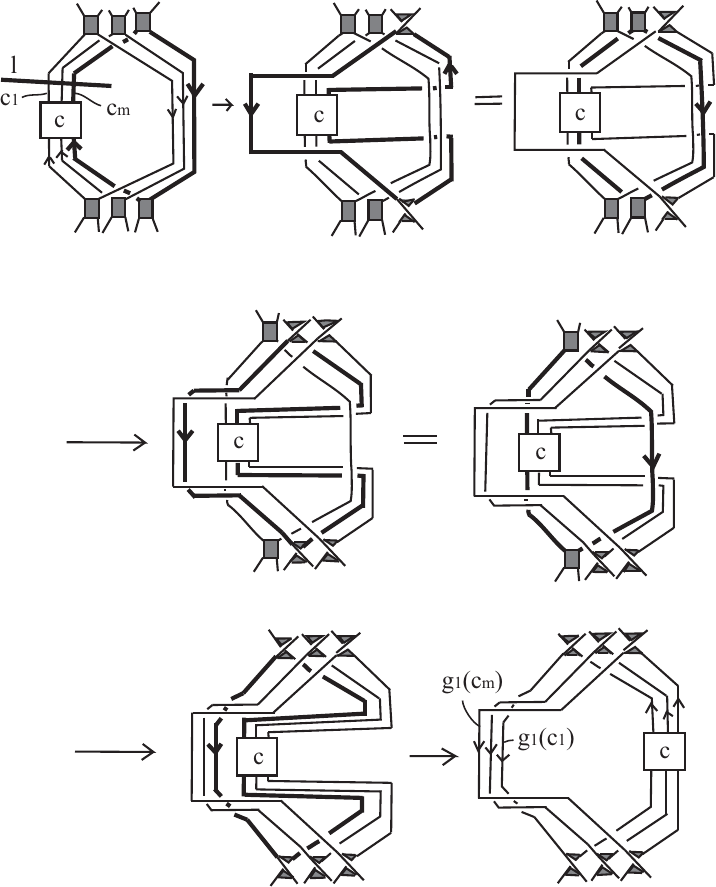}
  \caption{Deforming $S_{-1}$. }
  \label{0528-3}
\end{figure} 

\begin{figure} 
     \includegraphics*{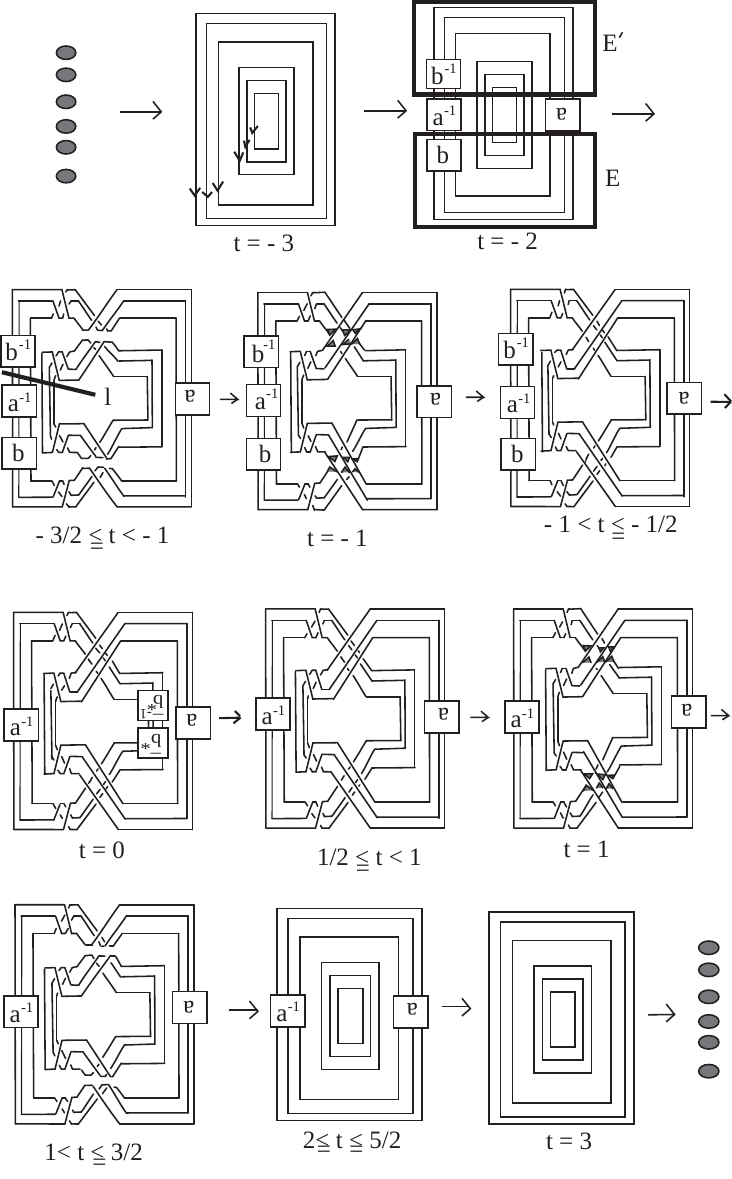}
  \caption{Motion picture of $S$ (5).}
  \label{B2-2}
\end{figure}   

\begin{figure} 
     \includegraphics*{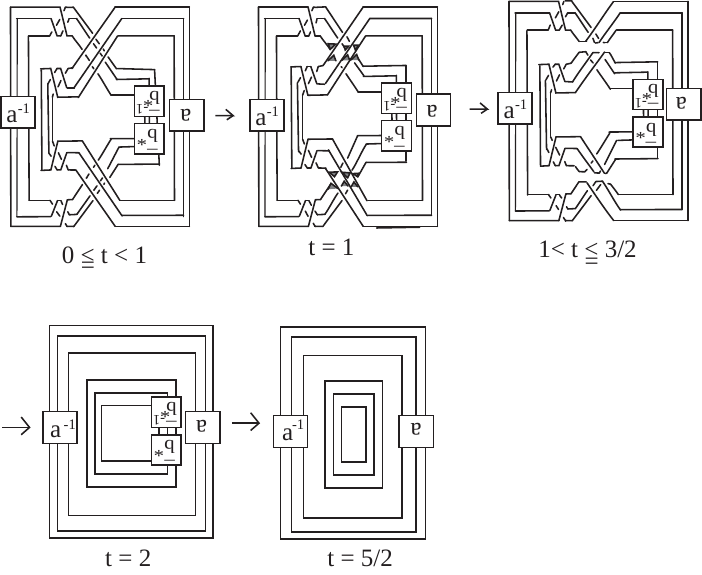}
  \caption{Motion picture of $S \cap \mathbb{R}^3[0, 5/2]$ (6).}
  \label{B2-3}
\end{figure}  

\begin{figure} 
     \includegraphics*{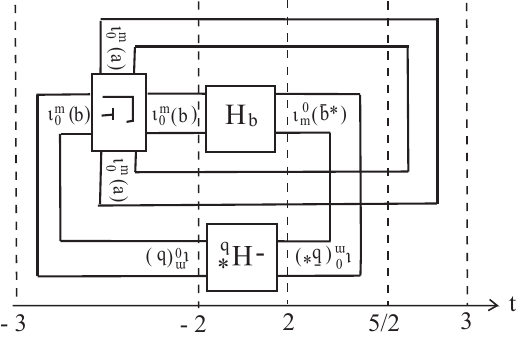}
  \caption{The $2m$-chart describing the surface braid $S_{\Gamma}$.}
  \label{B3-2}
\end{figure}   
{\bf Step 4}. \ 
We will show that $S$ is equivalent to the surface link whose 
 motion picture $\{S_t\}$ is as in 
Figure \ref{B2-2} for $t \in (- \infty, 0) \cup (5/2, +\infty)$, and 
otherwise 
as in Figure \ref{B2-3}. 
Further we will show that then $S$ is the closure of the required surface link $2m$-chart.  

Let us consider a trivial closed 1-braid $\mathrm{cl}(e_1)$ 
and two arcs in $\mathbb{R}^3$ 
with two untwisted bands connecting them as in Figure \ref{0528-2} (1). 
Let us take a 3-ball which contains $\mathrm{cl}(e_1)$ and the bands. 
Then, by an ambient isotopy of $\mathbb{R}^3$ 
relative the complement of the 3-ball, 
we can deform $\mathrm{cl}(e_1)$ and the bands to be as in Figure \ref{0528-2} (2). 

Consider a neighborhood of 
 $-\mathrm{cl}(\iota \bar{\Delta} \bar{\Delta}^{-1}) 
\cup \mathcal{B}_3$ 
in $S_{-1}$ of Step 3 (see Figure \ref{0613-1}). 
Since $-\mathrm{cl}(\iota \bar{\Delta} \bar{\Delta}^{-1})$ 
is 
equivalent to $-\mathrm{cl}(\iota e_m)$, 
we can take $m$ disjoint 3-balls in $\mathbb{R}^3$ 
such that the intersection of $S_{-1}$ and each 3-ball 
is as in Figure \ref{0528-2} (1). 
Let $c_j$ be the $j$th string of 
$-\mathrm{cl}(\iota \bar{\Delta} \bar{\Delta}^{-1})$. 
Apply the deformation of Figure \ref{0528-2} to the 3-balls which contains $\mathrm{cl}(c_1)$, 
$\mathrm{cl}(c_2), \ldots, \mathrm{cl}(c_m)$ 
and the attaching bands in this order (see Figure \ref{0528-3}, 
where we consider the case $c=e_m$).  
Let $\{g_u\}$ be the associated ambient isotopy. 
Since $c_{j-1}$ is over $c_{j}$ with respect to the $z$-axis, 
$g_1(c_{j})$ is over $g_1(c_{j-1})$ with respect to the $z$-axis ($j=2,3,\ldots, m$). 
Since $g_1(\iota\bar{\Delta})$ is the $m$-braid obtained from $\iota \bar{\Delta}$ by changing each crossing, 
for $\iota \bar{\Delta}$ with attaching bands, we have $g_1(\iota \bar{\Delta})=\iota \Delta$. 
Thus $g_1(\iota \bar{\Delta} \bar{\Delta}^{-1})=\iota \Delta \Delta^{-1}$. 
Further we can assume that 
$g_1(-\mathrm{cl}(\iota \bar{\Delta} \bar{\Delta}^{-1}))=
\mathrm{cl}(\iota \Delta^{-1} \Delta \Delta^{-1} \Delta)$ as in 
Figure \ref{0528-3}. 
Put 
$\mathcal{B}_4=g_1(\mathcal{B}_3)$, which is a band set such that each band has one twist. 
 Then $g_1(S_{-1})$ is a 
closed singular $2m$-braids such that the attaching band set $\mathcal{B}_4$ consists of 
$2m$ bands which 
correspond to the $m$ $\sigma_m$'s 
and $m$ $\sigma_m^{-1}$'s, connecting 
the inner closed $m$-braid and outer closed $m$-braid. 

Let $E_1$ be a cylinder in $\mathbb{R}^3$ indicated at $t=-1/2$ in Figure \ref{0613-1}. 
Note that $S_{-1/2} \cap E_1$ is 
the trivial $m$-braid. 
We can see that if we have an $m$-braid $c$ in $E_1$, then $g_1(c)=c$ in the cylinder 
$g_1(E_1)$ 
(see Figure \ref{0528-3}). 
Then by an ambient isotopy of (\ref{0528-9}) with $h_u=g_u$ and 
$t_1=-5$, $t_2=-4$, $t_3=4$, $t_4=5$, 
 we can deform $S$ to have the motion picture as in Figure \ref{B2-2}, as follows. 
  
 Let us denote by $e_{2m}$ the trivial $2m$-braid, and 
 let $\tilde{\mathcal{D}}_{-,2m}$ (resp. $\tilde{\mathcal{D}}_{+, 2m}$) be 
 a trivial disk system of $2m$ components in 
$\mathbb{R}^3(-\infty, t]$ (resp. $\mathbb{R}^3[t, +\infty)$) with boundary 
$\mathrm{cl}(e_{2m})[t]$, 
where $t \in \mathbb{R}$. 
 Let $\iota_0^m(\Gamma_T)$ be the $2m$-chart obtained from the $m$-chart 
 $\Gamma_T$ by adding $m$ trivial sheets after $\Gamma_T$, i.e. 
 the $2m$-chart presented by $\Gamma_T$. 
 And let $A(\iota_0^m(\Gamma_T))$ be 
 the simple braided surface over an annulus associated with $\iota_0^m(\Gamma_T)$. 
From now on, we denote by $\Delta$ (resp. $\Delta^\prime$) the $2m$-braid 
$\Delta_m$ (resp. $\Delta_m^\prime$). Further, we put $\Theta=\Theta_m$. 
Let us take $l$, which is indicated at $t=-3/2$ in Figure \ref{B2-2}. 
 Then the motion picture of $S$ is as follows (see Figure \ref{B2-2}), 
where we take the identified corresponding ends of a closed braid in 
$l$ for $t \in [-3,\, 3]$. 
 \begin{enumerate}
\item
$S \cap \mathbb{R}^3 (-\infty , -3]=\tilde{\mathcal{D}}_{-, 2m}$.  
We have 
$S_{-3}=
\mathrm{cl}(e_{2m})$. 

\item 
$S \cap \mathbb{R}^3 [-3, -2]=A(\iota_0^m(\Gamma_T))$. 
We have  
$S_{-2}=
\mathrm{cl}(\iota_0^m(a^{-1}bab^{-1}))$. 

 \item
 $S_{-2} \rightarrow S_{-3/2}$ is an isotopic transformation such that 
\[
S_{-3/2}=\mathrm{cl}( \iota^m_0 (a^{-1}b) \cdot \Delta^{\prime -1}
 \Delta^{-1} \Delta^{\prime} 
 \Delta  \cdot \iota^m_0 (a) \cdot \Delta^{-1} \Delta^{\prime -1} \Delta \Delta^{\prime} 
\cdot \iota^m_0 (b^{-1})). 
\]

 \item
$S_{-3/2} \rightarrow S_{-1/2}$ is a hyperbolic transformation along 
the band set $\mathcal{B}_4 \subset S_{-1}$. 
We have $S_{-1/2}=
\mathrm{cl}( \iota^m_0 (a^{-1} b)  \cdot \Delta^{\prime -1}
 \Delta^{-1} \Theta \cdot \iota^m_0 (a) \cdot \Theta^{-1} \Delta \Delta^{\prime} \cdot 
\iota^m_0 (b^{-1}))$. 

\item
$S_{-1/2} \rightarrow S_{0}$ is an isotopic transformation 
such that 
\[
S_{0}=\mathrm{cl}( \iota^m_0 (a^{-1}) \cdot \Delta^{\prime -1}  
 \Delta^{-1} 
\Theta \cdot \iota^m_0 (a) \cdot \iota_m^0(\bar{b}^* \bar{b}^{* -1}) 
\cdot \Theta^{-1} \Delta \Delta^{\prime}).
\]

\item
$S_{0} \rightarrow S_{1/2}$ is an isotopic transformation 
such that 
\[
S_{1/2}=\mathrm{cl}( \iota^m_0 (a^{-1}) \cdot \Delta^{\prime -1}  
 \Delta^{-1} 
\Theta \cdot \iota^m_0 (a) \cdot \Theta^{-1} \Delta \Delta^{\prime}).
\]
\item 
$S_{1/2} \rightarrow S_{3/2}$ is a hyperbolic transformation along the band set 
$\mathcal{B}_4 \subset S_1$. 
We have 
\[
S_{3/2}=\mathrm{cl}( \iota^m_0 (a^{-1}) \cdot \Delta^{\prime -1}  
 \Delta^{-1} \Delta^{\prime}  
 \Delta \cdot \iota^m_0 (a) \cdot \Delta^{-1}  
 \Delta^{\prime -1} \Delta \Delta^{\prime}).
\]

\item 
$S_{3/2} \rightarrow S_{2}$
is an isotopic transformation such that 
 $S_{2}=\mathrm{cl}(\iota_0^m(a^{-1}a))$. 
 
\item 
$S_t=S_2$ for $t \in [2, 5/2]$. 

 \item 
 $S_{5/2} \rightarrow S_{3}$ is an isotopic transformation such that 
  $S_3=\mathrm{cl}(e_{2m})$. 
 
\item 
$S \cap \mathbb{R}^3 [3, +\infty)=\tilde{\mathcal{D}}_{+,2m}$. 
\end{enumerate}

Let us consider an ambient isotopy $\{H_u\}_{u \in [0,1]}$ of $\mathbb{R}^3[0, 5/2]$ 
rel $\mathbb{R}^3[0] \cup \mathbb{R}^3[5/2]$ 
such that 
for each $u \in [0,1]$ and $\mathbf{x} \in \mathbb{R}^3$, 
\begin{equation} 
H_u(\mathbf{x}, t)=\begin{cases}
\big(h_{2ut}^{-1}(\mathbf{x}), t \big) & \mathrm{if} \ t \in [0, 1/2] \\
\big(h_{u}^{-1}(\mathbf{x}), t \big) & \mathrm{if} \ t \in [1/2, 2] \\
\big(h_{2u(t-2)} \circ h_u^{-1} (\mathbf{x}), t \big) & \mathrm{if} \ t \in [2, 5/2], 
\end{cases}
\end{equation}
where $\{h_u\}_{u \in [0,1]}$ is the ambient isotopy of $\mathbb{R}^3$ 
such that $h_u(S_0)=S_{u/2}$ for each $u$. 
By $\{H_u\}$, 
$S$ is equivalent to the surface link whose motion picture $\{S_t\}$ is 
as in Figure \ref{B2-2}  
for $t \in (-\infty, 0) \cup (5/2, \infty)$, and otherwise 
as follows (see Figure \ref{B2-3}).  
 \begin{enumerate} 
\item 
$S_{t}=S_0=\mathrm{cl}( \iota^m_0 (a^{-1}) \cdot \Delta^{\prime -1}  
 \Delta^{-1} \Theta \cdot 
\iota^m_0 (a) \cdot \iota^0_m (\bar{b}^* \bar{b}^{* -1}) \cdot \Theta^{-1} \Delta \Delta^{\prime})$ 
for $t \in [0,1/2]$. 

\item 
$S_{1/2} \rightarrow S_{3/2}$ is a hyperbolic transformation along the band set 
$\mathcal{B}_4 \subset S_1$. 
We have 
\[
S_{3/2}=\mathrm{cl}( \iota^m_0 (a^{-1}) \cdot \Delta^{\prime -1} 
 \Delta^{-1}  \Delta^{\prime} 
 \Delta \cdot \iota^m_0 (a) \cdot \iota^0_m (\bar{b}^* \bar{b}^{* -1}) \cdot \Delta^{-1} \Delta^{\prime -1} 
\Delta \Delta^{\prime}). 
\]

\item 
$S_{3/2} \rightarrow S_{2}$
is an isotopic transformation such that 
\[
 S_{2}=\mathrm{cl}(\iota_0^m(a^{-1}a) \cdot \iota_m^0(\bar{b}^* \bar{b}^{* -1})). 
\]
 
 \item 
 $S_2\rightarrow S_{5/2}$ is an isotopic transformation such that 
  $S_{5/2}=\mathrm{cl}(\iota_0^m(a^{-1} a))$. 
\end{enumerate}
Note that the band set $\mathcal{B}_4$ consists of $2m$ bands 
which correspond to the $m$ $\sigma_m$'s 
and $m$ $\sigma_m^{-1}$'s, connecting 
the inner and outer closed $m$-braids.  

The figures describing the motion pictures of $S$ are for the case $m=3$. 
However, since each step can be applied for an arbitrary positive integer $m$, 
we can regard them as describing the motion pictures of $S$ for 
arbitrary $m$. 
 Now, the motion picture $\{S_t\}$ of $S$ is as in 
Figure \ref{B2-2} for $t \in (- \infty, 0) \cup (5/2, +\infty)$, and 
otherwise 
as in Figure \ref{B2-3}. 
Thus $S$ is the closure of the surface braid 
 obtained from $S \cap \mathbb{R}^3[-3,3]$ by cutting it along $l [-3,3]$. 
 Let us denote the surface braid by $S_{\Gamma}$. 
Let us take cylinders $E$ and 
$E^{\prime}$ 
in $\mathbb{R}^3$ such that the closed $2m$-braid $S_t$ is presented by 
$\iota_0^m(a^{-1}) \cdot (S_t \cap E) \cdot \iota_0^m(a) \cdot (S_t \cap E^\prime)$ 
for $t \in [-2, 2]$, where $S_t \cap E$ and $S_t \cap E^\prime$ 
 are classical $2m$-braids; see the figure at 
$t=-2$ in Figure \ref{B2-2}. 
Then $S_{\Gamma} \cap (E [-2,2])$ is a braided surface of degree $2m$ 
such that 
\begin{enumerate}
\item 
$S_{-2} \cap E=\iota_0^m(b)$, 
\item
$S_{-2} \cap E \rightarrow S_{-3/2} \cap E$ 
is an isotopic transformation such that 
$S_{-3/2} \cap E =\iota_0^m(b) \cdot \Delta^{\prime -1} \Delta^{-1} \Delta^\prime \Delta$. 
\item 
$S_{-3/2} \cap E \rightarrow S_{-1/2} \cap E$ 
is a hyperbolic transformation along a band set 
which consists of $m$ bands corresponding to the $m$ $\sigma_m$'s. 
We have 
$S_{-1/2} \cap E =\iota_0^m(b) \cdot \Delta^{\prime -1} \Delta^{-1} \Theta$. 
\item
$S_{-1/2} \cap E \rightarrow S_{0} \cap E$ 
is an isotopic transformation such that 
$S_{0} \cap E =\Delta^{\prime -1} \Delta^{-1} \Theta \cdot \iota_m^0(\bar{b}^*)$. 
\item
$S_{0} \cap E \rightarrow S_{3/2} \cap E$ 
is a hyperbolic transformation along a band set 
which consists of $m$ bands corresponding to the $m$ $\sigma_m$'s. 
We have 
$S_{3/2} \cap E =\Delta^{\prime -1} \Delta^{-1} \Delta^\prime \Delta \cdot \iota_m^0(\bar{b}^*)$. 
\item
$S_{3/2} \cap E \rightarrow S_{2} \cap E$ 
is an isotopic transformation such that 
$S_{2} \cap E =\iota_m^0(\bar{b}^*)$. 
\end{enumerate}
Thus the braided surface $S_{\Gamma} \cap E[-2,2]$ is 
presented by the 1-handle $2m$-chart 
$H_b$. 
 If $a$ is the trivial $m$-braid, then $S \cap \mathbb{R}^3 [t]$ is equivalent to the 
trivial closed $2m$-braid for $t \in [-2, 2]-\{-1, 1\}$. 
Thus
$(S_t \cap E) \cdot (S_t \cap E^\prime)=e_{2m}$, 
and hence 
\begin{equation} \label{s}
S_t \cap E^\prime= (S_t \cap E)^{-1}
\end{equation}
 for $t \in [-2,2]-\{-1,1\}$. 
Further, $S_{-1} \cap E^\prime$ and $S_1 \cap E^\prime$ are singular classical $2m$-braids 
such that the singular points are presented by $m$-bands corresponding to 
$m$ $\sigma_m^{-1}$'s. 
Let us consider a $2m$-chart $\Gamma$. 
Let $\rho$ be a path which intersects with the edges of $\Gamma$ transversely. 
Let 
$c=\sigma_{i_1}^{\epsilon_{i_1}}\sigma_{i_2}^{\epsilon_{i_2}}\cdots \sigma_{i_\nu}^{\epsilon_{i_\nu}}$ 
be the $2m$-braid presented by $\rho \cap \Gamma$. 
Let us consider the orientation-reversed mirror image of $\Gamma$ and $\rho$, 
and denote it by $-\Gamma^*$ and $-\rho^*$ respectively. 
Then 
 $-\rho^* \cap (-\Gamma^*)$ presents 
$\sigma_{i_\nu}^{-\epsilon_{i_\nu}}\sigma_{i_{\nu-1}}^{-\epsilon_{i_{\nu-1}}}\cdots \sigma_{i_1}^{-\epsilon_{i_1}}$, 
which is $c^{-1}$. 
Hence, by (\ref{s}) we can see that $S_{\Gamma} \cap (E^{\prime} [-2,2])$ 
is the braided surface presented by a $2m$-chart $f(-H_b^*)$, 
where $-H_b^*$ is the orientation-reversed 
mirror image of $H_b$ and $f$ is the map which rotates $-H_b^*$ by $\pi$. 
More precisely, assuming that $-H_b^*$ is on a 2-disk $[-1,1] \times [-1,1]$,  
$f$ is a self-homeomorphism of the 2-disk 
defined by $f(t_1, t_2)=-(t_1, t_2)$. 
 Thus the surface braid $S_{\Gamma}$ 
  is presented by the surface link $2m$-chart as in Figure \ref{B3-2}. 
Since the $2m$-chart can be taken to be as in Figure \ref{0614-4} 
by an ambient isotopy, $S$ is presented by the required $2m$-chart. 
\qed 
\\

An orientable surface link $F$ is {\it trivial} (or {\it unknotted}) if 
there is an embedded 3-manifold $M$  with $\partial M=F$ 
such that each component of $M$ is a handlebody. 
An oriented surface link is called 
{\it ribbon} if it is obtained from a trivial 2-link $F_0$ by 1-handle surgeries 
along a finite number of mutually disjoint 1-handles attaching to $F_0$. 

 Kamada showed \cite{Kamada92} that surface links whose braid index is 
 at most three are indeed ribbon, and Shima showed \cite{Shima} that 
the turned spun $T^2$-knot of a non-trivial classical knot 
is not ribbon. Hence we obtain the following corollary.
\begin{cor} \label{c0614-7}
  Let $S$ be the torus-covering link presented by 
  a torus-covering $m$-chart. 
  Then the braid index of S is equal or less than $2m$. 
  In particular, the braid index of the turned spun $T^2$-knot 
  of the torus $(2,\,p)$-knot is four.
  \end{cor}

\begin{rem}
 Hasegawa \cite[Part 3 \lq\lq Chart description of twist-spun 
 surface-links"]{Hasegawa} showed that for the turned 
 spun $T^2$-link of a closed $m$-braid, 
 its braid index is at most $3m$.
\end{rem}
 
\section{Example} \label{0616-1}
\begin{figure}
    \includegraphics*{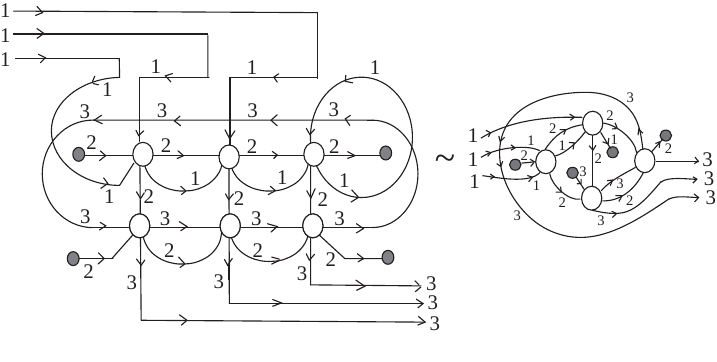}
\caption{The 1-handle $4$-chart $H_b$, where $b=\sigma_1^3$.}
\label{Fig3-11}
\end{figure}
 \begin{figure}
  \includegraphics*{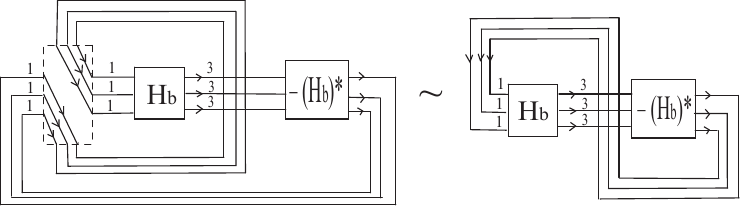}
\caption{The surface link $4$-chart $\Gamma_S$ obtained from $\Gamma_T$.}
\label{Fig3-12-3}
\end{figure}
 \begin{figure}
  \includegraphics*{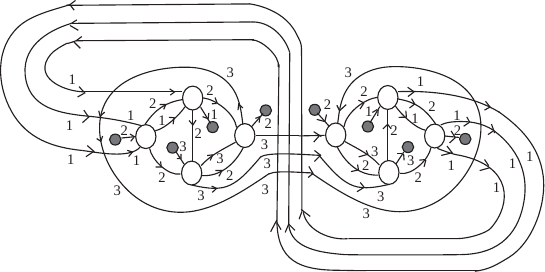}
\caption{The surface link $4$-chart $\Gamma_S$.}
\label{Fig3-12-2}
\end{figure}

Let us consider the torus-covering $2$-chart $\Gamma_T$ of Example \ref{0228-1} (2), 
that is, the torus-covering $2$-chart $\Gamma_T$ of Figure \ref{Fig3-10} (2). 
The torus-covering link $S$ associated with $\Gamma_T$ is 
 the turned spun $T^2$-knot of the right-handed trefoil. 
In the notations of Theorem \ref{theorem3-1}, we have
$a=b=\sigma_1^3$. 
By the theorem, the 1-handle $4$-chart $H_{b}$ is as in the left figure of Figure \ref{Fig3-11}. 
It is $C$-move equivalent to the right figure of Figure \ref{Fig3-11}. 
Thus the surface link $4$-chart 
$\Gamma_S$ on a 2-disk obtained 
from $\Gamma_T$ is as in the left figure of Figure \ref{Fig3-12-3}. 
It is equivalent to the right figure of Figure \ref{Fig3-12-3} 
by an ambient isotopy of the 2-disk. 
Thus $\Gamma_S$ is as in Figure \ref{Fig3-12-2}. 
\\

  \acknowledgement
The author would like to thank Professors Takashi Tsuboi and Elmar Vogt for suggesting this topic, 
and Professors Akio Kawauchi, Tomotada Ohtsuki, and the referee for their valuable advice.

\end{document}